\documentclass[3p,preprint]{elsarticle}

\usepackage{amsmath}
\usepackage{amsfonts}
\usepackage{amssymb,latexsym}
\usepackage{amsthm}
\usepackage{graphicx}

\newtheorem{theorem}{Theorem}[section]

\newtheorem{Remark}{Remark}[section]
\newtheorem{Corollary}{Corollary}[section]

\newtheorem{Proposition}{Proposition}[section]

\usepackage{color}
\usepackage{newlfont} 
\usepackage{subfigure}
\usepackage{verbatim}
\usepackage{rotating}
\usepackage{multirow}
\usepackage{units}
\usepackage{bm}
\usepackage[english]{babel}
\usepackage[utf8]{inputenc}
\usepackage{algorithm}
\usepackage[noend]{algpseudocode} 
\usepackage{booktabs}
\usepackage{caption}
\usepackage{hyperref}
\usepackage{subfigure}
\usepackage[export]{adjustbox}

\newcommand{\vol}{\mbox{vol}}

\journal{...}

\begin{document}

\begin{frontmatter}

\title{Effective numerical integration on complex shaped elements\\ by discrete signed measures}

\author[address-PD]{Laura Rinaldi}

\author[address-PD]{Alvise Sommariva}

\author[address-PD]{Marco Vianello\corref{corrauthor}}
\ead{marcov@math.unipd.it}

\cortext[corrauthor]{Corresponding author}

\address[address-PD]{University of Padova, Italy}

\begin{abstract}
We discuss a cheap and stable approach to polynomial moment-based compression of multivariate measures by discrete signed measures. The method is based on the availability of an orthonormal basis and a low-cardinality algebraic quadrature formula for an auxiliary measure in a bounding set. Differently from other approaches, no conditioning issue arises since no matrix factorization or inversion is needed. We provide bounds for the 
sum of the absolute values of the signed measure weights, 
and we make two examples: efficient quadrature on curved planar elements with spline boundary  (in view of the application to high-order
FEM/VEM), and compression of QMC integration on 3D elements with complex shape.
\end{abstract}

\begin{keyword}
MSC[2020]  65D32
\end{keyword}

\end{frontmatter}

\section{Introduction}
The possibility of compressing any finite measure into a low-cardinality discrete representing one, which has the some polynomial moments up to a certain degree, dates back to the famous Tchakaloff theorem, along with its successive developments; cf. e.g. \cite{T57,
Put97}. The implementation of Tchakaloff theorem has been an active research subject especially in the last decade, but the constraint of positivity of the weights makes quite costly the computational solutions via moment-matching, based on theoretical issues like Caratheodory theorem on conical combinations \cite{C11} or Davis-Wilhelmsen theorem on Tchakaloff sets \cite{W76}, together with suitable optimization techniques like linear or quadratic programming ; cf. e.g. \cite{KKN18,LL12,MP13,PSV17,RB15,SV15,T15} with the references therein. 

In the present paper we propose a general method still based on moment-matching (with a suitable orthogonal polynomial basis), that instead computes a {\em signed} discrete representing measure (that is, a quadrature formula with not all positive weights). We recently made a first step in this direction for integration on polyhedral elements \cite{SV25}, modifying and extending  the method proposed in \cite{LvPBD24}. Indeed, the present  approach has two remarkable features. First, the weights are computed simply (and cheaply) by a matrix-vector product followed by a vector scaling, so that {\em no matrix conditioning} issue can arise, even at high exactness degrees, because {\em no matrix inversion or factorization} is needed. Moreover, in some situations of practical interest, e.g. in the integration over elements for FEM/VEM, the relevant matrix can be computed once and for all, and {\em the weights, though not all positive, have a bounded 1-norm}, which is a necessary and sufficient condition for stability of the integration formula. Moreover, differently from \cite{SV25}, the stability analysis is here performed by rigorous estimates, based on a careful application of the classical Bessel, Cauchy-Schwarz and Jensen inequalities. 
As a result, full stability of the method is obtained, concerning both, the weights computation and the weights behavior.

The method is quite general, since either continuous or discrete measures can be compressed, as we show by the examples in the numerical section. The paper is organized as follows. The theoretical foundations and the relevant estimates are discussed in Section 2. In Section 3 we present two applications, namely 
cheap and stable quadrature on curved planar elements with spline boundary (in view of the application to high-order FEM/VEM), and compression of QMC integration on “difficult” 3D elements with complex  shape, obtained by basic set operations like intersection and union of simpler shapes. The corresponding Matlab and Python numerical codes are freely available.

\section{Moment-based discretized measure compression by orthogonal polynomials}

\begin{theorem}\label{thm2.1}
Let $\mu$ be a finite measure with support $\Omega\subseteq \mathbb{R}^d$, and 
$\lambda$ a finite measure with $\mathbb{P}_n$-determining support $B\subseteq \mathbb{R}^d$ (i.e., polynomials of total degree not exceeding $n$ which vanish there vanish everywhere in $\mathbb{R}^d$). Moreover, let $\{p_j\}_{1\leq j\leq N}$ be a $\lambda$-orthonormal basis for $\mathbb{P}_n$, $N=dim(\mathbb{P}_n)={n+d \choose d}$, 
and let $(X,\mathbf{u})=\{(x_i,u_i),\,1\leq i\leq M\}$, $M\geq N$, be the nodes $X$ and positive weights of a quadrature formula for $\lambda$, exact in $\mathbb{P}_{2n}$, that is 
\begin{equation} \label{exact2n}
\int_B{p(x)\,d\lambda}=\sum_{i=1}^M{u_i\,p(x_i)}\;,\;\;\forall p\in \mathbb{P}_{2n}\;.
\end{equation}
Denote by $V\in \mathbb{R}^{M\times N}$ the Vandermonde-like matrix $V=(v_{ij})=(p_j(x_i))$, by $D=D(\mathbf{u})=diag(u_i)$ the diagonal weight matrix, 
and by $\mathbf{m}=\{m_j\}\in \mathbb{R}^N$ the vector of moments
\begin{equation} \label{mom}
m_j=\int_\Omega{p_j(x)\,d\mu}\;,\;\;1\leq j\leq N\;.
\end{equation}

Then, the discrete signed measure $(X,\mathbf{w})=\{(x_i,w_i),\,1\leq i\leq M\}$ 
supported on $X\subseteq B$, with weights $\mathbf{w}=\{w_i\}\in \mathbb{R}^M$ computed as 
\begin{equation} \label{w}
\mathbf{w}=DV\mathbf{m}\;,
\end{equation}
gives an exact quadrature formula in $\mathbb{P}_n$ for $\mu$ with the following bound for the weights  
\begin{equation} \label{exact}
\int_\Omega{p(x)\,d\mu}=\sum_{i=1}^M{w_i\,p(x_i)}\;,\;\;\forall p\in \mathbb{P}_{n}\;,\;\;\|\mathbf{w}\|_1\leq \sqrt{\lambda(B)}\,\|\mathbf{m}\|_2\;.
\end{equation}
\end{theorem}
\vskip0.5cm
\noindent{\bf Proof.} First, we prove that the weight vector $\mathbf{w}$ is a solution of the underdetermined moment system $V^t \mathbf{w}=\mathbf{m}$, and thus the equality in (\ref{exact}) holds. Indeed, the matrix $D^{1/2}V$ is {\em orthogonal}, because in view of (\ref{exact2n}) we have 
$$
\sum_{i=1}^M{u_i\,p_h(x_i)\,p_k(x_i)}=\int_B{p_h(x)\,p_k(x)\,d\lambda}=\delta_{hk}\;,
$$
that is 
$(D^{1/2}V)^tD^{1/2}V=V^tDV=I$. Hence $V^t\mathbf{w}=V^tDV\mathbf{m}=\mathbf{m}$. 

On the other hand, $\|D^{-1/2}\mathbf{w}\|_2=\|D^{1/2}V\mathbf{m}\|_2=\sqrt{\mathbf{m}^tV^tDV\mathbf{m}}=\|\mathbf{m}\|_2$ and thus 
by {\em Cauchy-Schwarz inequality}
$$
\|\mathbf{w}\|_1=\sum_{i=1}^M{|w_i|}=\sum_{i=1}^M{(|w_i|/\sqrt{u_i}\;)\sqrt{u_i}}\leq \|\{|w_i|/\sqrt{u_i}\,\}\|_2\,\|\{\sqrt{u_i}\,\}\|_2
$$
$$
=\|D^{-1/2}\mathbf{w}\|_2\,\sqrt{\sum_{i=1}^M{u_i}}=\|\mathbf{m}\|_2\,\sqrt{\lambda(B)}\;.\hspace{1cm} \square
$$

\vskip0.5cm

\begin{Corollary}  
Let the assumptions of Theorem 2.1 be satisfied, with $\Omega\subseteq B$. 
If $B$ is bounded then 
\begin{equation} \label{BM}
\|\mathbf{w}\|_1\leq \sqrt{\lambda(B)}\,\mu(\Omega)\,\sqrt{\max_{x\in  B}K_n(x,x)}\;\;,
\end{equation}
where $K_n(x,x)=\sum_{j=1}^N{p_j^2(x)}$ is the so called Christoffel polynomial or also reciprocal Christoffel function of $\lambda$ (the diagonal of its reproducing kernel). 
On the other hand, if  $\mu$ and $\lambda$ are absolutely continuous (w.r.t. the Lebesgue measure) with densities $\omega\in L^1(\Omega)$ and $\sigma\in L^1(B)$, i.e. $d\mu=\omega(x)\,dx$ and $d\lambda=\sigma(x)\,dx$, and $\omega^2/\sigma\in L^1(\Omega)$, then 
\begin{equation} \label{stab}
\|\mathbf{w}\|_1\leq \sqrt{\lambda(B)}\,\sqrt{\|\omega^2/\sigma\|_{L^1(\Omega)}}\;\;.
\end{equation}
\end{Corollary}
\vskip0.5cm
\noindent{\bf Proof.} 
As for the first bound, in view of (\ref{exact}) it is sufficient to observe that by {\em Jensen inequality} for probability measures
$$
\left(\int_\Omega{p_j(x)\frac{d\mu}{\mu(\Omega)}}\right)^2\leq \int_\Omega{p_j^2(x)\frac{d\mu}{\mu(\Omega)}}
$$
and thus
$$
\|\mathbf{m}\|_2^2=\sum_{j=1}^N{\left(\int_\Omega{p_j(x)\,d\mu}\right)^2}\leq \mu(\Omega)\sum_{j=1}^N{\int_\Omega{p_j^2(x)\,d\mu}}=\mu(\Omega)\int_\Omega{K_n(x,x)\,d\mu}
$$
$$
\leq \mu^2(\Omega)\max_{x\in \Omega}K_n(x,x)
\leq\mu^2(\Omega)\max_{x\in B}K_n(x,x)\;.
$$
Observe that the maxima exist since $B$ and $\Omega$ are closed (being supports of measures) and bounded, and hence compact. 

Concerning the second bound, denoting by $I_\Omega(x)$ the indicator function of $\Omega\subseteq B$, we can write 
$$
m_j=\int_\Omega{p_j(x)\,d\mu}=\int_\Omega{p_j(x)\,\omega(x)\,dx}=\int_B{p_j(x)\,I_\Omega(x)\,\frac{\omega(x)}{\sigma(x)}\,\sigma(x)\,dx}
$$
which shows that the $\{m_j\}$ are the coefficients of the Fourier expansion in the orthogonal polynomial basis of the function {\textcolor{black}{$I_\Omega\,\omega/\sigma$, that belongs to $L^2_\lambda(B)$}} since 
$$
\int_B{\left(I_\Omega(x)\,\frac{\omega(x)}{\sigma(x)}\right)^2\sigma(x)\,dx}=\int_B{I_\Omega^2(x)\,\frac{\omega^2(x)}{\sigma(x)}\,dx}=
\int_\Omega{\frac{\omega^2(x)}{\sigma(x)}\,dx}\;.
$$
Consequently, by {\em Bessel inequality}
$$
\|\mathbf{m}\|_2^2=\sum_{j=1}^N{m_j^2}\leq 
\|I_\Omega\,\omega/\sigma\|^2_{L^2_\lambda(B)}=\|\omega^2/\sigma\|_{L^1(\Omega)}\;.
\hspace{1cm} \square$$

\vskip0.5cm

Some remarks are now in order, to deepen the properties of the measure compression method studied above. 

\begin{Remark}
The upper bound in (\ref{stab}) is independent of $n$, thus ensuring stability of the quadrature formula in the usual sense. On the other hand, whatever is $\mu$ the right-hand side of (\ref{BM}) has sub-exponential growth in $n$ (a weak form of quadrature stability) if $\lambda$ is a Bernstein-Markov measure, that holds for instance when $\lambda$ is absolutely continuous; cf. e.g. \cite{BLPW15}. To make an example, if $\lambda$ is the product Chebyshev measure on $B=[-1,1]^2$, it is known that $\sqrt{\max_{x\in [-1,1]^2}K_n(x,x)}\leq \frac{1}{\pi}\,\sqrt{2n^2+2n+1}\lesssim 0.45\,n$, 
cf. \cite{DMSV14}.
\end{Remark}
\vskip0.5cm
\begin{Remark}
Though the previous results have a simple proof, resting on basic linear algebra and inequalities, they are quite general and give some deep consequences. Indeed, by the signed measure $(X,\mathbf{w})$ we can ``compress'' either continuous or discrete measures $\mu$, 
provided that their polynomial moments are computable. 

A first unusual consequence is for example that we can integrate exactly polynomials somewhere, by sampling ``elsewhere''. Indeed, the sets $\Omega$ and $B$ could be even disjoint; a problem is however that we expect the polynomials $\{p_j(x)\}$ and $K_n(x,x)$, and thus also $\|\mathbf{m}\|_2$, to increase rapidly with $n$ when 
$\Omega\cap (\mathbb{R}^d\setminus B)\neq \emptyset$; cf. e.g. \cite{BMW20,L24,LPP22} with the references therein. Another (well-known) fact is that also moment computation can be performed in some situations by sampling only on the boundary $\partial\Omega$, by the divergence theorem, as we shall see in the application to quadrature on curved planar elements. So that we can integrate a polynomial on a domain with few or even no samples inside the domain and one might wonder to which extent this can be extended by the present theory to other functions, well approximated by polynomials  (a partial answer is given by Proposition 2.1 below). In practice however, as we have seen in the absolutely continuous case, a sufficient condition for quadrature stability, i.e. for $\|\mathbf{w}\|_1$ having a bound independent of $n$, is that  
$\Omega\subseteq B$, even if the weights are not all positive. 

Another consequence of Theorem 2.1 is that the weights are computed simply (and cheaply) by a matrix-vector product, so that no matrix conditioning issue can arise, even at high exactness degrees. This makes the present approach quite different from other moment-based methods for the compression of discrete measures, that need Vandermonde matrices factorization or inversion and/or numerical optimization techniques, especially those implementing Tchakaloff theorem; see, e.g.,  \cite{KKN18,LvPBD24,LL12,PSV17,RB15,SV15,T15} with the references therein. 
\end{Remark}
\vskip0.5cm

The measure compression method developed above allows to integrate exactly polynomials in $\mathbb{P}_n$. But in which terms could it be applied to the integration of continuous functions? Which is the cumulated effect of errors in the computation of moments and errors in the sampled function values? We treat both these aspects in the case $\Omega\subseteq B$, by resorting to the general result proved in \cite[\S2]{SVZ08} on the role of moment errors in algebraic quadrature. Indeed, denote by $\tilde{\mathbf{m}}\approx \mathbf{m}$ the approximate moments, by $\tilde{\mathbf{w}}=DV\tilde{\mathbf{m}}$ the resulting approximate weights, by $E_n(f;B)=\inf_{p\in \mathbb{P}_n}{\|f-p\|_{\infty,B}}$ the best uniform approximation error of $f\in C(B)$, by $\mathbf{f}=\{f(x_i)\}$ the exact sampled values of $f$ and by $\tilde{\mathbf{f}}\approx \mathbf{f}$ the perturbed sampled values. Moreover, let $\langle\cdot\,,\cdot\rangle$ denote the scalar product 
in $\mathbb{R}^M$. Then 
we have the following
\begin{Proposition}
Let the assumptions of Theorem 2.1 be satisfied with $\Omega\subseteq B$ and $B$ bounded. 

Then the following cumulative error estimate holds 
$$
\left|\int_\Omega{f(x)\,d\mu}-\langle\tilde{\mathbf{w}},\tilde{\mathbf{f}}\rangle\right|\leq \left(\mu(\Omega)+\|\mathbf{w}\|_1\right)\,E_n(f;B)+\|\mathbf{w}\|_1\,\|\mathbf{f}-\tilde{\mathbf{f}}\|_\infty
$$
\begin{equation} \label{errest}
+\left(\|f\|_{L^2_\lambda(B)}+\sqrt{\lambda(B)}\,(E_n(f;B)+\|\mathbf{f}-\tilde{\mathbf{f}}\|_\infty)\right)\,\|\mathbf{m}-\tilde{\mathbf{m}}\|_2\;,
\end{equation}
where $\|\mathbf{w}\|_1$ can be bounded as in (\ref{BM}), or as in (\ref{stab} when both $\mu$ and $\lambda$ are absolutely continuous.

\end{Proposition}
\vskip0.5cm
\noindent{\bf Proof.} 
Using formula (15) in \cite[\S2]{SVZ08} we get
$$
\left|\int_\Omega{f(x)\,d\mu}-\langle\tilde{\mathbf{w}},\mathbf{f}\rangle\right|\leq \left(\mu(\Omega)+\|\mathbf{w}\|_1\right)\,E_n(f;B)+\left(\|f\|_{L^2_\lambda(B)}+\sqrt{\lambda(B)}\,E_n(f;B)\right)\,\|\mathbf{m}-\tilde{\mathbf{m}}\|_2\;.
$$
On the other hand,
$$
\left|\int_\Omega{f(x)\,d\mu}\langle\tilde{\mathbf{w}},\tilde{\mathbf{f}}\rangle\right|\leq \left|\int_\Omega{f(x)\,d\mu}-\langle\tilde{\mathbf{w}},\mathbf{f}\rangle\right|+\left|\langle\mathbf{w},\mathbf{f}-\tilde{\mathbf{f}}\rangle\right|+\left|\langle\mathbf{w}-\tilde{\mathbf{w}},\mathbf{f}-\tilde{\mathbf{f}}\rangle\right|
$$
$$
\leq  \left|\int_\Omega{f(x)\,d\mu}-\langle\tilde{\mathbf{w}},\mathbf{f}\rangle\right|+\|\mathbf{w}\|_1\,\|\mathbf{f}-\tilde{\mathbf{f}}\|_\infty+\|\mathbf{w}-\tilde{\mathbf{w}}\|_1\,\|\mathbf{f}-\tilde{\mathbf{f}}\|_\infty\;,
$$
where reasoning as in the proof of Theorem 2.1 we get
$\|\mathbf{w}-\tilde{\mathbf{w}}\|_1\leq \sqrt{\lambda(B)}\,
\|\mathbf{m}-\tilde{\mathbf{m}}\|_2$. 
Putting together all the bounds we obtain eventually estimate (\ref{errest}). $\hspace{0.5cm} \square$
\vskip0.5cm
\begin{Remark}
When $\|\mathbf{w}\|_1$ is bounded, estimate (\ref{errest}) shows that the quadrature formula $\langle\mathbf{w},\mathbf{f}\rangle$ converges to the integral for any continuous function (the condition being indeed also necessary, by the  multivariate extension of Polya-Steklov theorem). 
We recall that the rate of $E_n(f;B)$ can be estimated depending on the regularity of $f$, if $B$ is a so-called ``Jackson
compact set'', i.e. a Jackson like inequality is available on $B$. This means that for every $\ell>0$ there exists a positive integer
$k_\ell$ and a constant $C_\ell(f)$, depending on the partial derivatives of $f$ up to order $k_\ell$, such that for $f\in C^{k_\ell}(B)$ the estimate  
$E_n(f;B)\leq C_\ell(f)\,n^{-\ell}$
holds for $n>\ell$. Examples are balls with $k_\ell=\ell$ and boxes with $k_\ell=\ell+1$; cf. \cite{P09}.   

Notice that in practice, by neclecting the products of errors (which tendentially give a much smaller contribute) , we get the approximate bound 
\begin{equation} \label{dom} 
\left|\int_\Omega{f(x)\,d\mu}-\langle\tilde{\mathbf{w}},\tilde{\mathbf{f}}\rangle\right|\lesssim \left(\mu(\Omega)+\|\mathbf{w}\|_1\right)\,E_n(f;B)+\|\mathbf{w}\|_1\,\|\mathbf{f}-\tilde{\mathbf{f}}\|_\infty
+\|f\|_{L^2_\lambda(B)}\,\|\mathbf{m}-\tilde{\mathbf{m}}\|_2\;.
\end{equation}
\end{Remark}

\section{Numerical examples}

In this numerical section we consider the application of the {\it{Cheap}} technique to the computation of integrals on planar spline curvilinear elements, potentially useful for solving PDEs problems by high-order FEM of VEM, and to the compression of discrete measures with large cardinality support such as those corresponding to QMC integration on
3D elements with complex shape.
{\color{black}In each example we will 
\begin{itemize}
\item discuss how to compute the moments;
\item check numerically the ADE ({\em Algebraic Degree of Exactness});
\item display the average cputime for determining the {\sl{Cheap}} rules, comparing them with other methods.
\end{itemize}
}
The open source codes that accomplish the numerical experiments of this section are available on GitHub (respectively at {\cite{AS}} for the Matlab version, and at {\cite{LRS}} for the Python one).

\subsection{Cheap and stable quadrature on planar spline curvilinear elements}

As first case we consider the application of {\it Cheap} technique to planar spline curvilinear elements.
These are Jordan domains $\Omega \subset {\mathbb{R}}^2$ whose boundary $\partial \Omega$ is described by parametric equations $$x=\phi(t)\;,\;y=\psi(t)\;,\;t \in [a,b]\;,\;\phi, \psi \in C([a,b])\;,\;\phi(a)=\phi(b)\;,\; 
	\psi(a)=\psi(b)\;,$$
    and there is a partition $\{I^{(k)}\}_{k=1,\ldots,M}$ of $[a,b]$ and partitions $\{I_j^{(k)}\}_{j=1,\ldots,m_k}$ of each $I^{(k)}$, such that the restrictions of  $\phi, \psi$ on each $I^{(k)}$ are splines of degree $\delta_k$, w.r.t. the subintervals $\{I_j^{(k)}\}_{j=1,\ldots,m_k}$.
 {\color{black}Namely, given the arc ``vertices''} $V_k \in \partial \Omega$, $k=1,\ldots,M+1$, then $\partial {{\Omega}}:=\cup_{k=1}^M V_k \frown V_{k+1}$, each curved side $V_k \frown V_{k+1}$ is tracked by a spline curve of degree $\delta_k$, interpolating an ordered subsequence of {\color{black}control} knots $P_{1,k}=V_k,P_{2,k},\ldots,P_{m_k-1,k},P_{m_k,k}=V_{k+1}$.

Elements of this kind have been considered in \cite{SV09}, where we introduced a panelization technique to determine a rule with ADE equal to $n$, of PO type, that is with positive weights and some nodes possibly outside the element. Such quadrature rules have been often used in the recent FEM and VEM literature, in view of their flexibility, but have the drawback of producing a large number of nodes per element. 

Later on, in \cite{SV21} we proposed a moment-based algorithm that determines PI rules, i.e. with positive weights and all nodes in the element. To this purpose, we have developed a fast {\em in-domain} routine, to allocate points from which low-cardinality rules are extracted. The drawback here is the necessity of computing a sparse nonnegative solution to the underdetermined moment-matching system, where we used the Lawson-Hanson algorithm as implemented in the Matlab {\tt lsqnonneg} function, or suitable accelerated variants based on the concept of ``deviation maximization'' instead of column pivoting for the underlying QR factorizations; cf. \cite{DODM23,DMV20}. In any case, this part turns out to be the computational bulk of the method, and leads to a significant cost per element.

In the {\it Cheap} implementation, we proceed differently, avoiding at all matrix factorizations and system solution, as well as sub-tessellation of the element. 
From this point of view the method is also alternative to the triangulation-free approach recently proposed in \cite{LvPBD24}, {\color{black}and gives a first extension to curved elements of the approach recently adopted in \cite{SV25} for linear polyhedra}. Instead, we compute:
\begin{itemize}
\item the smallest cartesian rectangle $B$ containing the element $\Omega$;
\item the tensorial Gauss-Chebyshev rule ($\mbox{ADE}=2n$) on the rectangle with $M=(n+1)^2$ nodes $X\subset B$;
\item the modified Chebyshev moments 
\begin{equation} \label{mom2d}
m_j=\int_\Omega{p_j(x,y)\,dxdy}=\oint_{\partial\Omega}{P_j(x,y)\;dy}\;,\;\;\frac{\partial P_j}{\partial x}=p_j\;,\;\;j=1,\ldots,{\color{black}N=\frac{(n+1)(n+2)}{2}}\;,
\end{equation}
by Gauss-Green theorem and Gauss-Legendre quadrature on the polynomial sub-arcs of $\partial \Omega$, where $p_j(x,y)=\tau_{i_1}(x)\tau_{i_2}(y)$,  
$0\leq i_1+i_2\leq n$, is the suitably 
(e.g. lexicographically) ordered orthonormal Chebyshev basis of $B$;
\item the weights {\bf{w}} of the cheap rule, relatively to the nodes $X$, by a single matrix-by-vector product as in (\ref{w}).
\end{itemize}

\begin{figure}[h]
	\centering
	\hspace{-1cm}
    
	\includegraphics[scale=0.465,clip,valign=t]{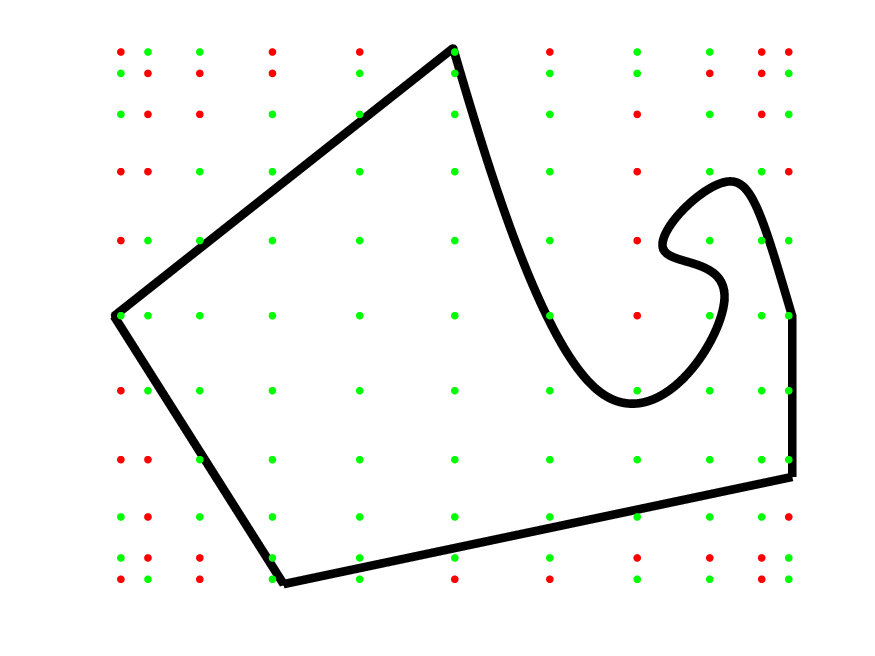}
	\hspace{0.2cm}
	\includegraphics[scale=0.56,clip,valign=t]{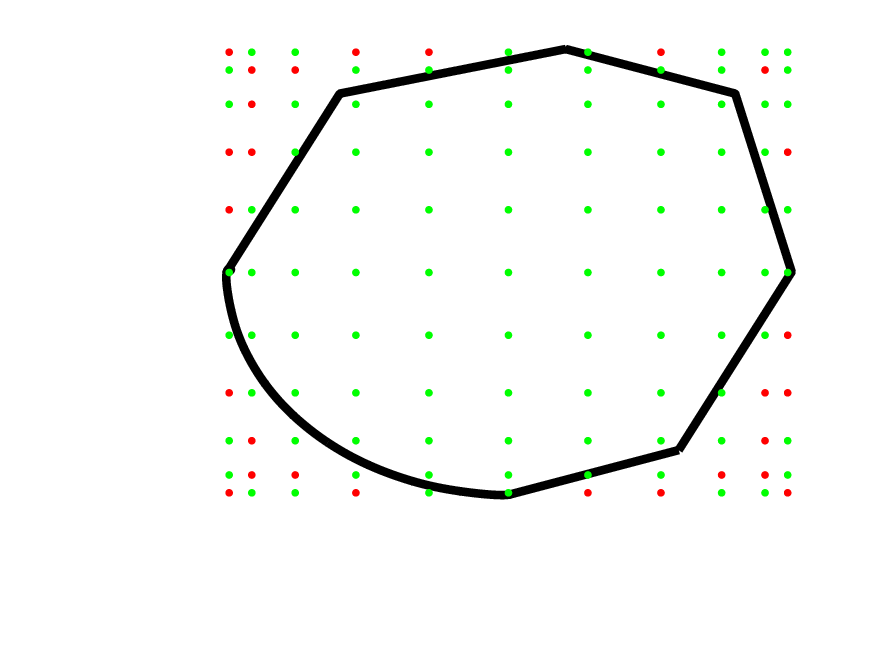}
	\caption{The planar curvilinear elements $\Omega_1$ and $\Omega_2$ and the nodes of a cheap formula with $\mbox{ADE}=10$. Green dots: nodes with positive weights; red dots: nodes with negative weights.}
	\label{fig:01A}
\end{figure}

In the examples below we test the Matlab routines on a computer with a $2.7$ GHz Intel Core $i5$ CPU, with $16$ GB of RAM.
In particular, we consider the curvilinear elements $\Omega_1$ and $\Omega_2$ depicted in {\ref{fig:01A}}, where 
\begin{itemize}
\item $\Omega_1$ {\color{black} is nonconvex} with four linear sides (linear splines) and {\color{black} one} curved side ({\color{black}cubic spline});
\item $\Omega_2$ {\color{black} is convex} with six linear sides (linear splines) and {\color{black} one} curved side ({\color{black}cubic spline}). 
\end{itemize}
Elements of this kind could typically appear when intersecting elements of a polygonal mesh with the domain near a curved boundary {\color{black}(tracked by splines)}.

\begin{figure}[h]
	\centering
	\hspace{-1cm}
	\begin{minipage}{0.35\linewidth}
    \centering
	\includegraphics[scale=0.45,clip]{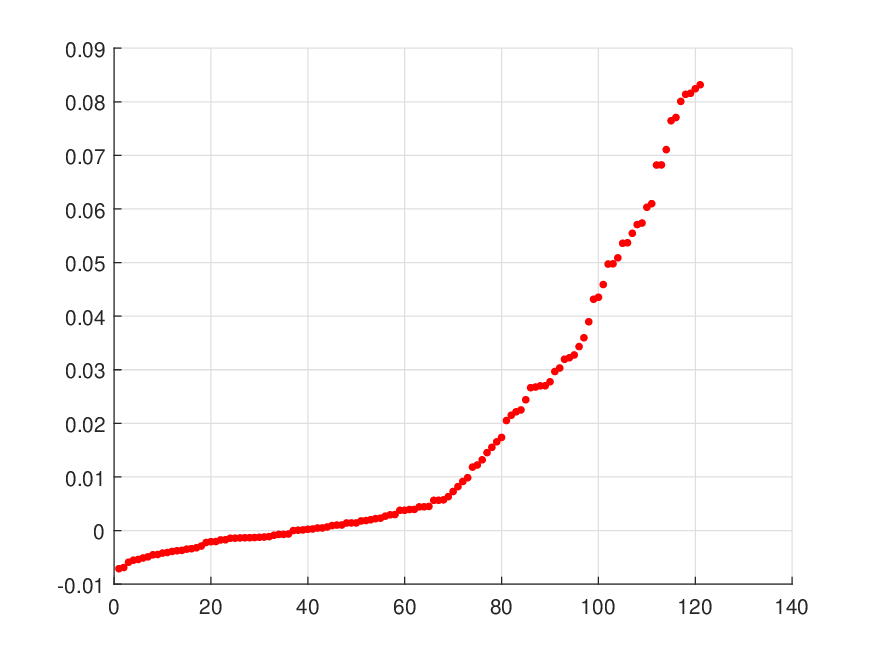}
     \end{minipage}
	\hspace{0.2cm}
	\begin{minipage}{0.4\linewidth}
		\centering
	\includegraphics[scale=0.45,clip]{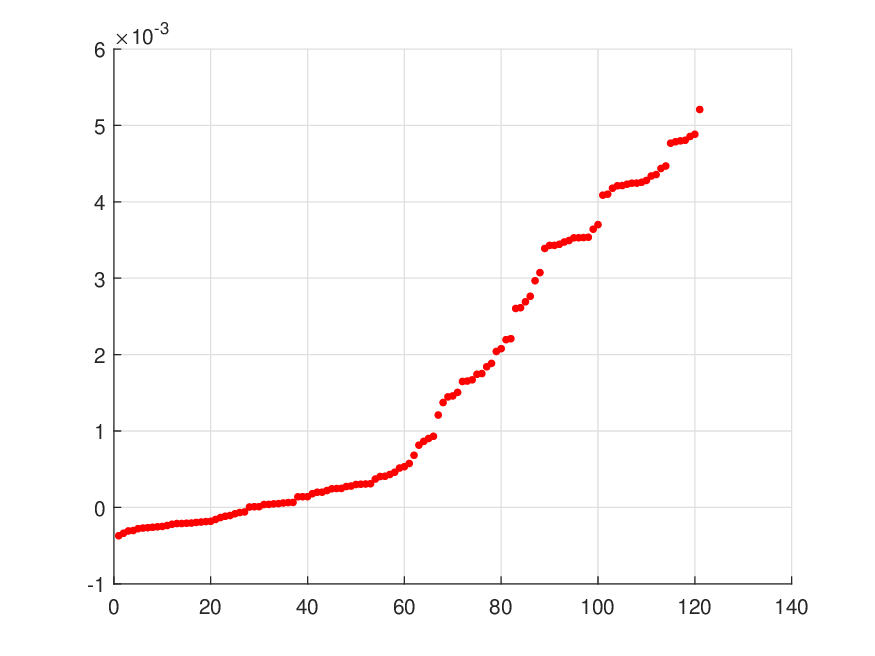}
\end{minipage}
	\caption{$121=11^2$ weights of the cheap rule with $\mbox{ADE}=10$ on the planar curvilinear elements $\Omega_1$ (left) and $\Omega_2$ (right), in increasing order.}
	\label{fig:01B}
\end{figure}

In Figure {\ref{fig:01A}}, we plot the nodes of the cheap rules with $\mbox{ADE}=10$ equal to 10 (green dots express positive weights, while red dots non positive ones). In Figure {\ref{fig:01B}} we report the corresponding weights in increasing order.

\begin{figure}[h]
	\centering
	\hspace{-1cm}
	\begin{minipage}{0.35\linewidth}
		\centering
		\includegraphics[scale=0.45,clip]{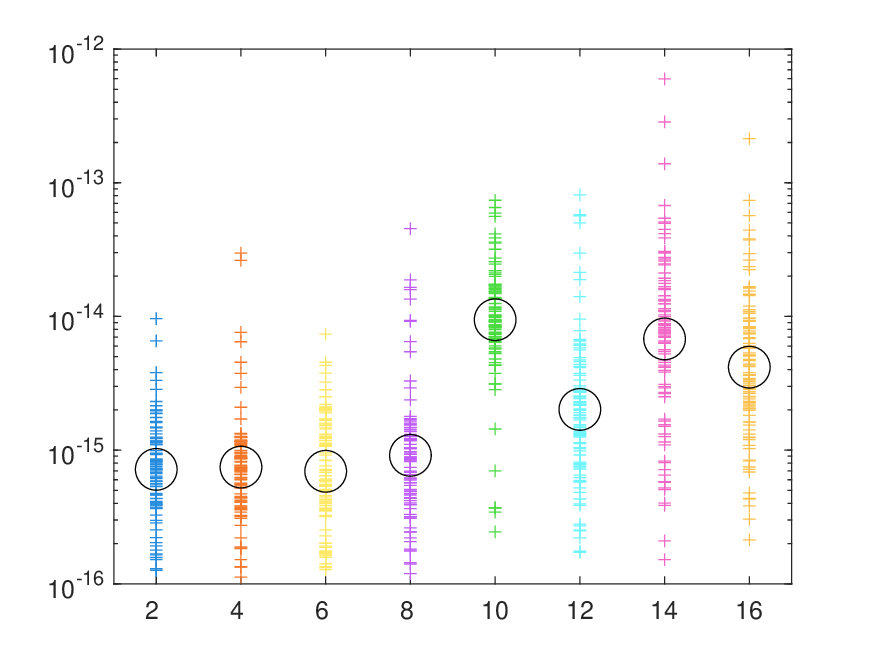}
	\end{minipage}
	\hspace{0.5cm}
	\begin{minipage}{0.4\linewidth}
		\centering
		\includegraphics[scale=0.45,clip]{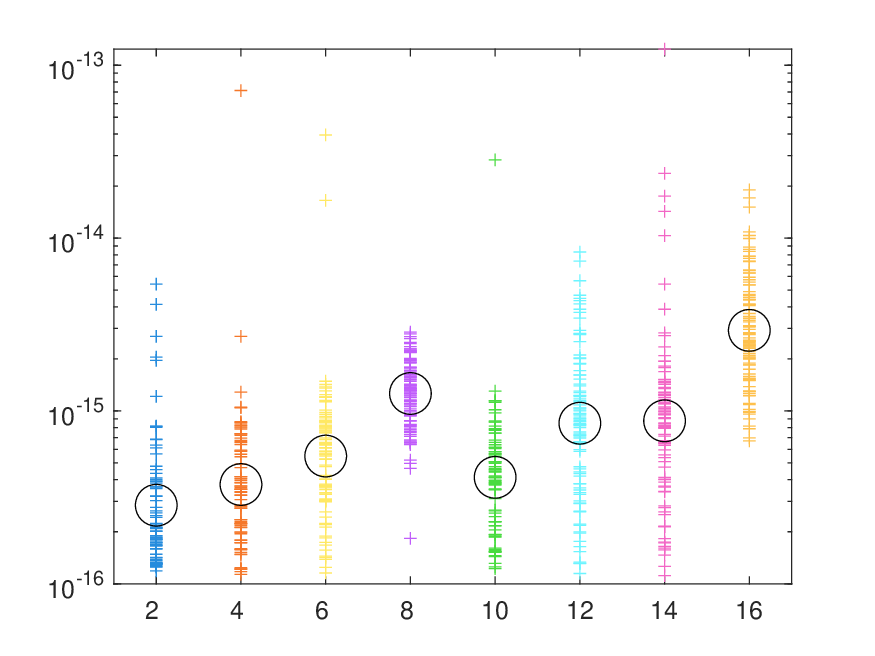}
	\end{minipage}
	\caption{{\color{black}Small crosses: relative quadrature errors for 100 trials of random polynomials $(c_0+c_1 x+c_2 y)^n$ on the elements $\Omega_1$ (left) and $\Omega_2$ (right). Circles: geometric mean of the relative errors. The abscissae are the ADE of the formulas.}}
	\label{fig:02}
\end{figure}

Next, to check the exactness of these rules, for any fixed algebraic degree of exactness $n=2,4,\,\ldots,16$ we compute {\color{black}100 trials of the integrals}
$$
I_{s,n}=\int_{\Omega_s} (c_0+c_1 x +c_2 y)^n\,dx\,dy\;,\;\;s=1,2\;
$$
with uniform random coefficients $c_j \in(0,1)$, $j=0,1,2$. The reference values of the integrals have been computed by applying Gauss-Green theorem and Gauss-Legendre quadrature as in {\cite{SV21}}.
In Figure \ref{fig:02} we have plotted by small crosses the relative errors made in the trials and by a circle their geometric mean. The results show that the implemented rules have mean errors not far from machine precision.
 
Next we take into account the quality of numerical integration of the functions 

\begin{eqnarray}
f_1(x,y)&=& \exp(-(x^2+y^2)),\\ 
f_2(x,y)&=& (x^2+y^2)^{11/2},\\ 
f_3(x,y)&=& (x^2+y^2)^{3/2},\\ 	
\end{eqnarray}

\noindent that present different regularity: $f_1$ is analytic and entire, whereas $f_2$ and $f_3$ have singularities of certain derivatives at the point $(0,0)$ in $\Omega_1$ and 
$\Omega_2$. In Table \ref{tab:01} we list the relative errors of the cheap rules, varying the ADE. {\color{black} The reference values of the integrals have been computed by the routine with ADE equal to 30.}

\begin{table}
	\centering
\begin{tabular}{| c || c | c | c ||| c | c | c ||}
	\hline
	 ADE  & $f_1$ & $f_2$ & $f_3$ & $f_1$ & $f_2$ & $f_3$\\
	 \hline
 	   2  & 3e-02 & 6e-01 & 2e-02 & 7e-05 & 3e-01 & 4e-02 \\
 	   4  & 5e-03 & 2e-01 & 3e-04 &  9e-07 & 6e-01 & 3e-03 \\
 	   6  & 1e-04 & 4e-03 & 2e-05 &  2e-09 & 2e-02 & 2e-05 \\
 	   8  & 5e-05 & 2e-04 & 5e-07 &  4e-12 & 2e-03 & 8e-05 \\
 	  10  & 4e-06 & 5e-07 & 2e-07 &  1e-14 & 5e-05 & 5e-05 \\
 	  12  & 2e-08 & 1e-09 & 7e-08 &  2e-15 & 5e-06 & 2e-06 \\
 	  14  & 7e-09 & 1e-10 & 3e-08 &  9e-16 & 5e-07 & 1e-05 \\
 	  16  & 5e-10 & 1e-11 & 7e-09 &  6e-16 & 5e-08 & 2e-06 \\

\hline
\end{tabular}
\caption{{\color{black}Relative errors of the cheap rules for the integration of $f_1,f_2,f_3$ on the elements $\Omega_1$ (left) and $\Omega_2$ (right).}}
\label{tab:01}
\end{table}

We conclude the section by displaying the {\color{black} construction cputime and the cardinality of the rules, for the two examples above; see Tables 2 and 3.} As anticipated, we compare the cheap rules with specialized rules for elements with spline boundary, namely rules of PO type introduced in \cite{SV09} and of PI type described in \cite{SV21}.
Concerning the cputime, we run 100 tests over $\Omega_1$ and $\Omega_2$. The startup process for cheap algorithm (construction of the bivariate Chebyshev-Vandermonde matrix) has been skip since it can be done once and for all for any mesh.
The results show that the cputime order of the present procedure is comparable to that of \cite{SV09}, with the advantage of a much lower cardinality of the rule. {\color{black}The PI rules proposed in \cite{SV21} have lower cardinality, but their computation is more expensive, especially at higher ADE.}

Finally, we report in Table {\ref{TAB_SPL_W}} the stability parameters $\|\mathbf{w}\|_1/|\sum_j w_j|\approx \|\mathbf{w}\|_1/vol(\Omega_s)$, $s=1,2$. The experiments show that these quantities are close to the optimal value $1$, so confirming quantitatively the expected stability of the rule.

\begin{Remark} 
Such a combination of no conditioning problem, low construction cost, stability and low cardinality, makes the Cheap rules a promising alternative for the time-consuming integral computations of FEM/VEM discretizations with curved polygonal elements, in the construction of stiffness and mass matrices.
\end{Remark}

\begin{table}
	\centering
	\begin{tabular}{| c || c | c | c ||| c | c | c ||}
		\hline 
		
		ADE  & CH & PO & PI & CH & PO & PI \\
		\hline
 	   2  & 5e-04 & 9e-04 & 9e-04 & 7e-04 & 9e-04 & 2e-03 \\
 	   4  & 4e-04 & 8e-04 & 2e-03 & 6e-04 & 8e-04 & 2e-03 \\
 	   6  & 4e-04 & 8e-04 & 7e-03 & 6e-04 & 8e-04 & 3e-03 \\
 	   8  & 4e-04 & 8e-04 & 6e-03 & 6e-04 &  8e-04 & 6e-03\\
 	  10  & 4e-04 & 8e-04 & 1e-02 & 6e-04 &  8e-04 & 1e-02 \\
 	  12  & 5e-04 & 8e-04 & 3e-02 & 7e-04 & 8e-04 & 3e-02 \\
 	  14  & 5e-04 & 8e-04 & 5e-02 & 9e-04 &  8e-04 & 5e-02 \\
 	  16  & 6e-04 & 8e-04 & 8e-02 & 8e-04 &  8e-04 & 7e-02 \\

\hline
\end{tabular}
\caption{{\color{black}Comparison of average construction cputimes (in seconds) of the cheap rule (CH), with the PO rule from \cite{SV09} and the PI rule from \cite{SV21}, on the elements $\Omega_1$ (left)  and $\Omega_2$ (right).}}
\label{tab:02}
\end{table}

\begin{table}
	\centering
	\begin{tabular}{| c || c | c | c || c | c | c ||}
		\hline 
		
		ADE  & CH & PO & PI & CH & PO & PI\\
		\hline
 	   2  &       9  &     76  &       6  &       9  &     132  &       6  \\
 	   4  &      25  &     171  &      15  &     25  &    297  &      15  \\
 	   6  &      49  &    304  &      28  &     49  &    528  &      28  \\
 	   8  &      81  &    475  &      45  &      81  &    825  &      45  \\
 	  10  &     121  &    684  &      66  &     121  &    1188  &      66  \\
 	  12  &     169  &    931  &      91  &     169  &    1617  &      91  \\
 	  14  &     225  &    1216  &     120  &     225  &    2112  &     120  \\
 	  16  &     289  &    1539  &     153  &    289  &    2673  &     153  \\

\hline
\end{tabular}
\caption{{\color{black}As in Table 2 concerning the cardinality of the rules.}}
\label{tab:03}
\end{table}

\begin{table}
	\centering
	\begin{tabular}{| c || c | c | c || c | c | c ||}
		\hline 
		
		ADE  & CH($\Omega_1$) & CH($\Omega_2$) \\
		\hline
 	   2  & 1.15 & 1.05 \\
 	   4  & 1.14 & 1.05 \\
 	   6  & 1.15 & 1.08 \\
 	   8  & 1.12 & 1.05 \\
 	  10  &  1.10 & 1.06  \\
 	  12  &  1.10 & 1.04 \\
 	  14  &  1.09 & 1.05 \\
 	  16  &  1.09 & 1.04 \\
\hline
\end{tabular}
\caption{The stability parameters $\|\mathbf{w}\|_1/|\sum_j w_j|$ of the cheap rules on the elements $\Omega_1$ (left)  and $\Omega_2$ (right).}
\label{TAB_SPL_W}
\end{table}

\subsection{{\color{black} Compression of discrete measures with large support: QMC on complex 3D shapes}}

In this section we consider the compression of discrete measures with large support on 3D geometries that may have a complex shape. The {\color{black} context} is the following. Given {\color{black}a sequence of low-discrepancy points} $\{Q_i\}_{i=1,\ldots,L}$ on {\color{black}a compact domain} $\Omega \subset {\mathbb{R}}^3$, one may apply the Quasi-Montecarlo (QMC) method to approximate the integral {\color{black}w.r.t. the Lebesgue measure} of a continuous function $f$, that is
$$
{\color{black}\int_\Omega f(x,y,z)\,dx\,dy\,dz \approx \frac{ {\vol(\Omega)}}{L}\, \sum_{i=1}^L f(Q_i)\;}.
$$
Typically one generates a {\color{black}low-discrepancy} sequence of cardinality $K$ on a parallelepiped $B$ that contains 
$\Omega$, and then by an in-domain routine {\color{black}(implementing the indicator function of $\Omega$)} extracts those $L$ points {\color{black}$\{Q_i=(x_i,y_i,z_i)\}_{i=1,\ldots,L}$} belonging to $\Omega$.
If the {\color{black}volume} of $\Omega$ is not easily available, as it happens for domains with complex geometrical shape, since
${\color{black}\vol(\Omega) \approx \vol(B) L/K}$, 
one can consider the approximation
\begin{equation} \label{QMC}
{\color{black}\int_\Omega{f(x,y,z)\,dx\,dy\,dz}\approx \int_\Omega{f(x,y,z)\,d\mu}{\textcolor{black}{:=}}\frac{\vol(B)}{K}\,\sum_{i=1}^L f(x_i,y_i,z_i)\;.}
\end{equation}
that is, QMC integration can be seen as computing the integral of $f$ w.r.t. a discrete measure $\mu$ defined by the nodes $Q_i$ and the equal weights ${\color{black}\nu_i=\vol(B)/K}$, $i=1,\ldots,L$. At this point, as described in Theorem {\ref{thm2.1}}, one can easily achieve {\color{black}a cheap rule with $M$ points in $B$}, where in general $M \ll L$, that is a compression of the original QMC rule. {\color{black}The resulting cheap QMC rule is exact for $f\in \mathbb{P}_n$, and integrates $f\in C(\Omega)$ with an additional error term of the order of the best polynomial approximation in 
$\mathbb{P}_n$ to $f$ on $\Omega$. We refer the reader e.g. to \cite{DP10} on QMC integration theory and to \cite{ESV22,ESV24} for a discussion on the error of compressed QMC (notice that a different and more costly compression method is there adopted, based on moment-matching via Nonnegative Least-Squares).} 

To give the idea, we consider two elements $\Omega_3$ and $\Omega_4$ with a complex shape. In particular 
\begin{itemize}
\item $\Omega_3$ is the intesection of the unit ball with a nonconvex polyhedron with 20 vertices and {\textcolor{black}{30 triangular facets}};
\item $\Omega_4$ is the union of $5$ balls ${\cal{B}}(C_k,r)$ {\color{black}with centers $\{C_k\}_{k=1,\ldots,5}$ in $[0,1]^3$ and equal radius $r=0.5$}.
\end{itemize}

\begin{figure}[h]
	\centering
	\hspace{-1cm}
    
\includegraphics[width=0.33\textwidth,valign=c]{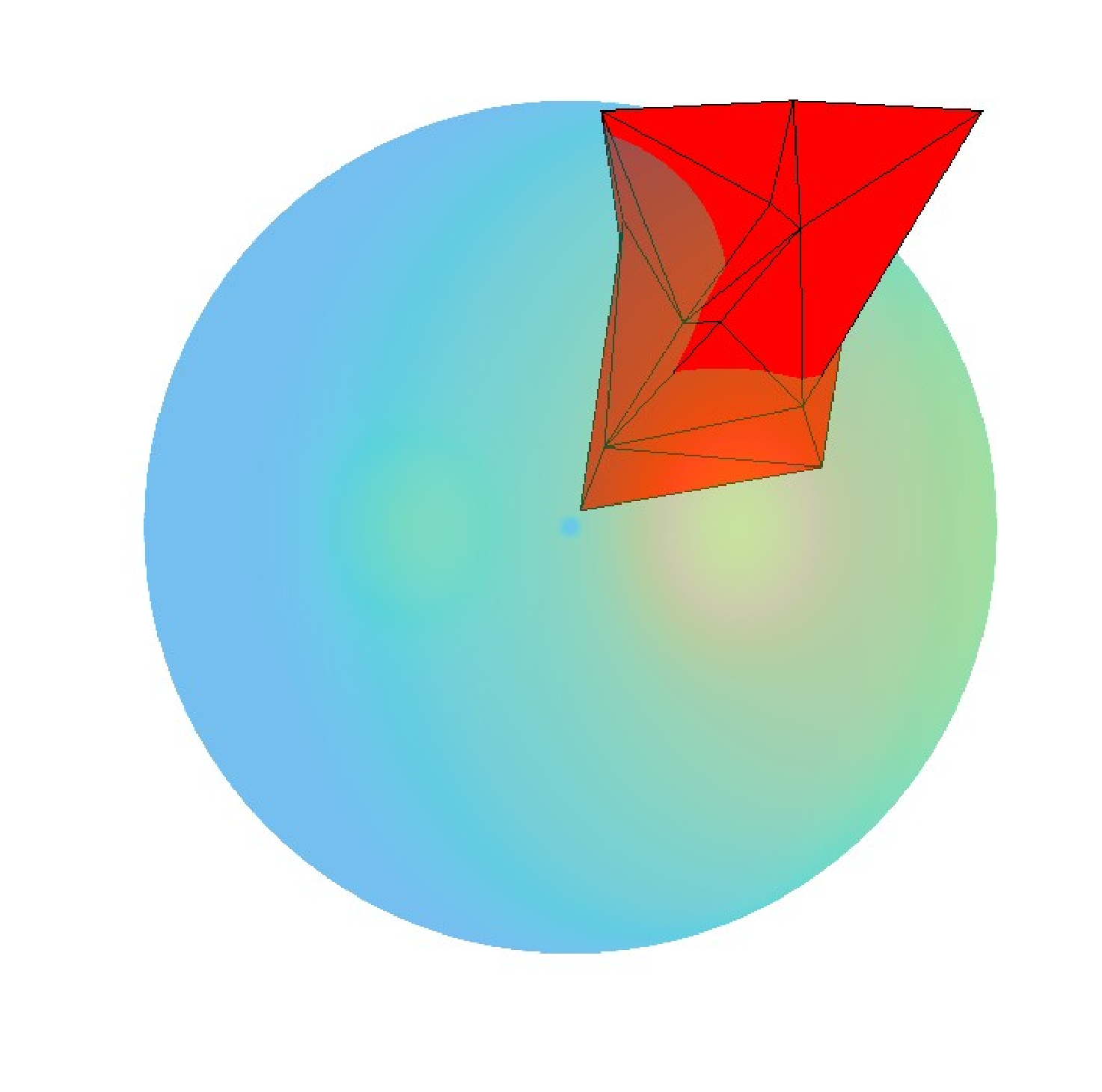}
\hspace{2cm}
\includegraphics[width=0.24\textwidth,valign=c]{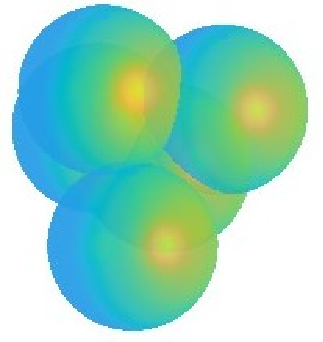}
	\caption{The multivariate elements $\Omega_3$ and $\Omega_4$.}
	\label{fig4}
\end{figure}

{\color{black}As it can be seen in Figure {\ref{fig4}}, the shapes are particularly complex and it is not easy to determine algebraic quadrature rules on these elements, since they would require an accurate tracking of the specific element geometry. This difficulty makes QMC integration an appealing alternative. To construct a QMC rule 
\begin{itemize}
\item for $\Omega_3$, we set as bounding box $B$ the parallelepiped obtained by intersection of cube $[-1,1]^3$ (circumscribed to the sphere), with the smaller parallelepiped containg the polyhedron (easily available from the vertices); 
\item for $\Omega_4$, we set as bounding box $B$ the smaller parallelepiped that contains all the five circumscribed cubes to the balls defining the element.
\end{itemize}

Next, to determine a cheap QMC rule we compute 
\begin{itemize}
\item the tensorial Gauss-Chebyshev rule ($\mbox{ADE}=2n$) on the parallelepiped with $M=(n+1)^3$ nodes $X\subset B$;
\item the QMC Chebyshev moments 
\begin{equation} \label{QMCmom)}
m_j=\int_\Omega{p_j(x,y,z)\,d\mu_{qmc}}\,:=\frac{\vol(B)}{K}\,\sum_{i=1}^L p_j(x_i,y_i,z_i)\;,\;\;j=1,\ldots,{\color{black}N=\frac{(n+1)(n+2)(n+3)}{6}}\;,
\end{equation}
where $p_j(x,y,z)=\tau_{i_1}(x)\tau_{i_2}(y)\tau_{i_3}(z)$,  
$0\leq i_1+i_2+i_3\leq n$, is the suitably 
(e.g. lexicographically) ordered orthonormal Chebyshev basis of $B$;
\item the weights {\bf{w}} of the cheap QMC rule, relatively to the nodes $X$, by a single matrix-by-vector product as in (\ref{w}).
\end{itemize}

In our numerical experiments we first scale $K=10^5$ Halton points in the bounding box $B$. Since in-domain functions are trivially available for balls and have been implemented by efficient algorithms for general polyhedra (cf. e.g. \cite{H15}), we can determine the QMC rules for both the elements. Such formulas have respectively $L=23076$ and $L=37379$ nodes.

To test the quality of QMC compression by the cheap rule, we plot the relative errors on 100 trials of the QMC rules applied to random polynomials of degree $n=2,4,\ldots,16$, with uniform random coefficients $c_j \in (0,1)$
$$
I_{s,n}=\int_{\Omega_s}{(c_0+c_1 x+c_2 y + c_3 z)^n\,d\mu_{qmc}} 
=\frac{\vol(B)}{K}\,\sum_{i=1}^L(c_0+c_1 x_i+c_2 y_i + c_3 z_i)^n\;,\;\;s=3,4.
$$
}

\begin{figure}[h]
	\centering
	\hspace{-1cm}
	\begin{minipage}{0.35\linewidth}
    \centering
	\includegraphics[scale=0.45,clip]{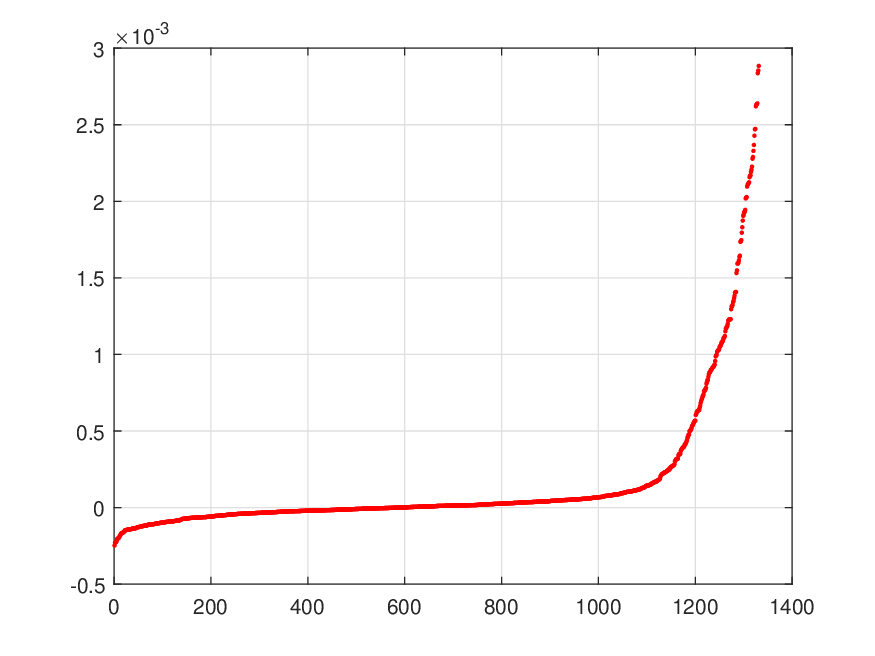}
     \end{minipage}
	\hspace{0.2cm}
	\begin{minipage}{0.4\linewidth}
		\centering
	\includegraphics[scale=0.45,clip]{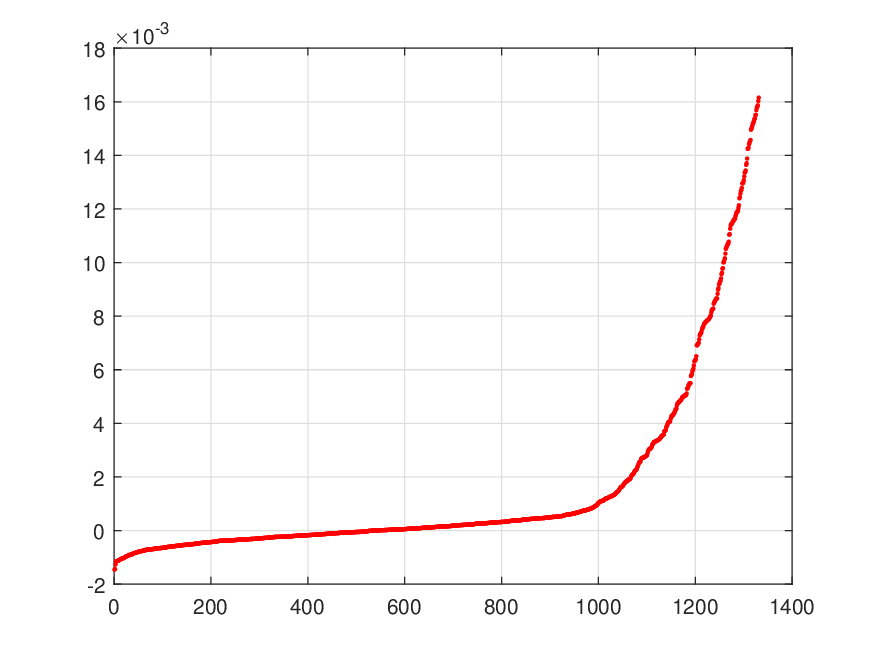}
\end{minipage}
	\caption{Weights in increasing order of the cheap QMC rules with $\mbox{ADE}=10$ on the 3D elements $\Omega_3$ (left) and $\Omega_4$ (right). {\color{black}The corresponding nodes are  $11^3=1331$, much less than the hundreds of thousands of points in the original QMC rules.}}
	\label{fig:04B}
\end{figure}

\begin{figure}[h]
	\centering
	\hspace{-1cm}
	\begin{minipage}{0.35\linewidth}
		\centering
		\includegraphics[scale=0.45,clip]{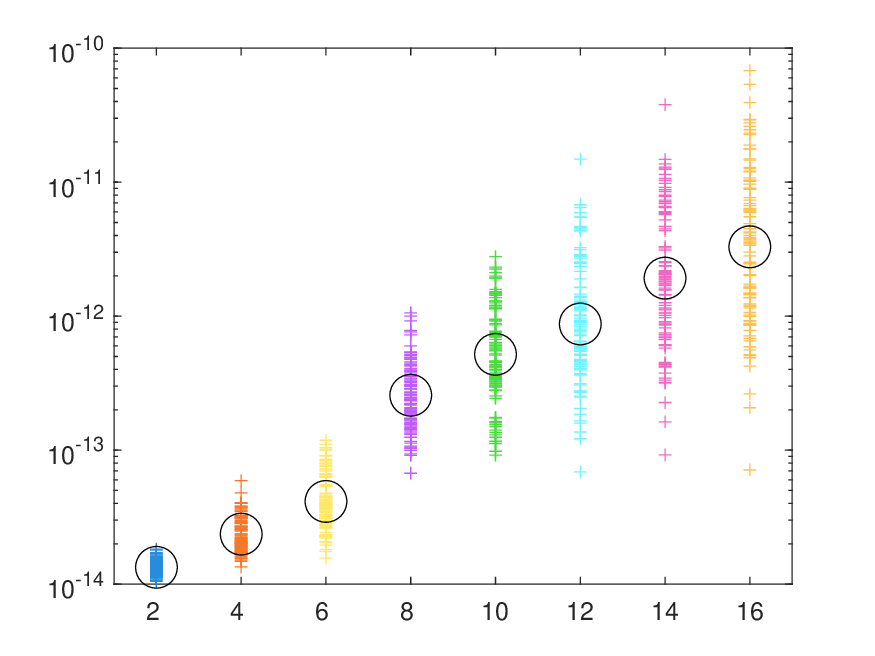}
	\end{minipage}
	\hspace{0.5cm}
	\begin{minipage}{0.4\linewidth}
		\centering
		\includegraphics[scale=0.45,clip]{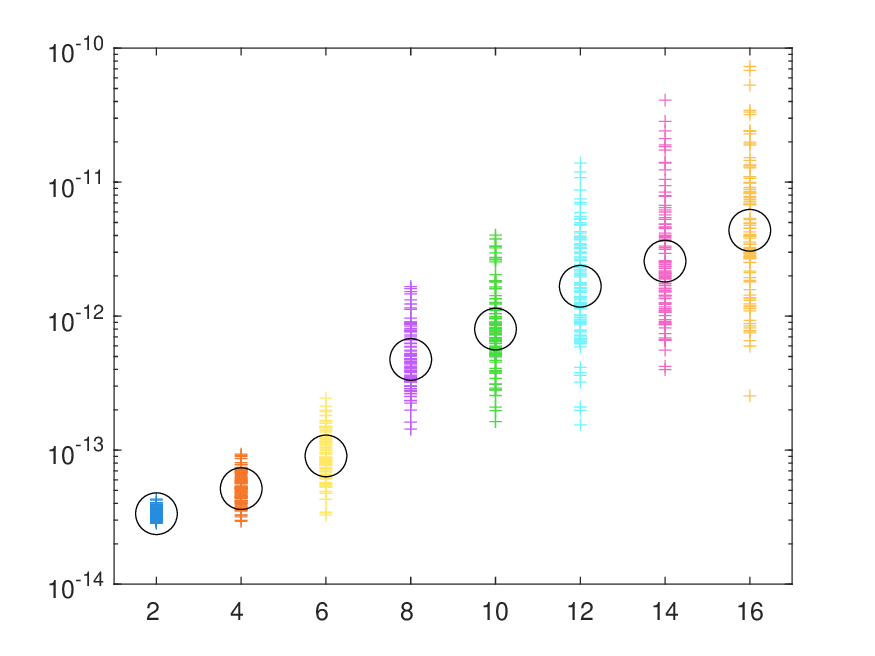}
	\end{minipage}
	\caption{{\color{black}Small crosses: relative quadrature errors for 100 trials of random polynomials $(c_0+c_1 x+c_2 y+c_3 z)^n$ on the elements $\Omega_3$ (left) and $\Omega_4$ (right). Circles: geometric mean of the relative errors. The abscissae are the ADE of the formulas.}}
	\label{fig:05}
\end{figure}

In Figure \ref{fig:05} we have plotted by
small crosses the relative errors made in the trials and by a circle their geometric mean. The results show that the implemented rules have mean errors that range from $10^{-14}$ to $10^{-12}$ depending on $n$. 

Next, in Table \ref{tab:04} we examine the quality of numerical integration of the functions 

\begin{eqnarray}
f_1(x,y,z)&=& \exp(-(x^2+y^2+z^2)),\\ 
f_2(x,y,z)&=& ((x-x_0)^2+(y-y_0)^2+(z-z_0)^2)^{11/2},\\ 
f_3(x,y,z)&=& ((x-x_0)^2+(y-y_0)^2+(z-z_0)^2)^{3/2}, 	
\end{eqnarray} 
with $(x_0,y_0,z_0)=(0.51,0.26,0.63)\in \Omega_3$ or   $(x_0,y_0,z_0)=(0.21,0.36,0.51)\in \Omega_4$. Again, the test functions present different regularity: $f_1$ is analytic and entire, whereas $f_2$ and $f_3$ have singularities of certain derivatives at $(x_0,y_0,z_0)$.  As reference values of the integrals, we considered those obtained by QMC rules on $\Omega_3$ and $\Omega_4$, constructed starting with $10^6$ points in the respective bounding boxes.

\begin{table}
\centering
\begin{tabular}{| c || c | c | c ||| c | c | c || c |}
\hline
	 ADE  & $f_1$ & $f_2$ & $f_3$ & $f_1$ & $f_2$ & $f_3$ & {\mbox{card}}\\
\hline
 2   &  3e-02   &  8e+00   &  1e-01  &  6e-02   &  5e+01   &  2e-01 & 27 \\
 4   &  1e-03   &  2e+00   &  6e-03  &  6e-03   &  2e+01   &  2e-02 & 125 \\
 6   &  5e-05   &  6e-02   &  5e-04  &  7e-04   &  4e-01   &  1e-03 & 343\\
 8   &  4e-05   &  1e-02   &  2e-04  &  7e-04   &  4e-01   &  2e-03 & 729\\
 10  &  4e-05   &  6e-03   &  3e-04  &  6e-04   &  7e-05   &  2e-03 & 1331\\
 12  &  4e-05   &  5e-03   &  4e-04  &  6e-04   &  9e-03   &  2e-03 & 2197\\
 14  &  4e-05   &  5e-03   &  4e-04  &  6e-04   &  9e-03   &  2e-03 & 3375 \\
 16  &  4e-05   &  5e-03   &  4e-04  &  6e-04   &  9e-03   &  2e-03 & 4913\\
\hline
\hline
QMC & 4e-05   &  5e-03   &  4e-04 &  6e-04   &  9e-03   &  2e-03  & \\
\hline
\end{tabular}
\caption{Relative errors of the cheap QMC rules for the integration of $f_1,f_2,f_3$ on the elements $\Omega_3$ (left) and $\Omega_4$ (right). The last row shows the relative errors of the original QMC rules with $10^5$ points in the bounding box $B$. In the last column we list the cardinality of the cheap rules that must be compared with the $L=23076$ and $L=37379$ nodes of the QMC rules on $\Omega_3$ and $\Omega_4$, respectively.}
\label{tab:04}
\end{table}

\begin{table}
	\centering
	\begin{tabular}{| c || c | c |}
		\hline 
		
		ADE \, & CH($\Omega_3$) & CH($\Omega_4$) \\
		\hline
 	   2  & 3e-03 & 5e-03 \\
 	   4  & 9e-03 & 1e-02 \\
 	   6  & 2e-02 & 3e-02 \\
 	   8  & 3e-02 & 5e-02 \\
 	  10  &  6e-02 &  9e-02  \\
 	  12  & 9e-02 & 1e-01 \\
 	  14  & 1e-01 & 2e-01 \\
 	  16  & 2e-01 & 3e-01 \\
            \hline
 {\mbox{QMC}}  & 1e+00 & 1e-01\; \\
 \hline
            
\end{tabular}

\caption{Average construction cputimes of the cheap QMC rule (in seconds), on the elements $\Omega_3$ (left)  and $\Omega_4$ (right). The last row shows the average cputimes for determining the nodes of the original QMC rules.}
\label{tab:05B}
\end{table}

\begin{table} 
\centering
\begin{tabular}{| c || c | c | c || c | c | c ||}
		\hline 
		ADE  & CH($\Omega_3$) & CH($\Omega_4$) \\
		\hline
 	   2  & 1.54 & 1.52 \\
 	   4  & 1.76 & 1.45 \\
 	   6  & 1.57 & 1.28 \\
 	   8  & 1.36 & 1.28 \\
 	  10  &  1.30 & 1.21  \\
 	  12  &  1.29 & 1.18 \\
 	  14  &  1.26 & 1.19 \\
 	  16  &  1.23 & 1.16 \\
\hline
\end{tabular}
\caption{The stability parameter $\|\mathbf{w}\|_1/|\sum_j w_j|$ of the cheap rules on the elements $\Omega_3$ (left)  and $\Omega_4$ (right).}
\label{TAB_QMC_W}
\end{table}

In Table \ref{tab:05B} we list the construction cputime, for the two examples above. As anticipated, these times do not include the construction of the tensorial Gauss-Chebyshev rule that can be done once, independently of the element.
After having at hand the cheap QMC rules, we run 100 tests over $\Omega_3$ and $\Omega_4$. 
The average cputimes of the cheap rules on the two elements are quite similar, while those of QMC are not, taking respectively $1$ and $0.1$ seconds. We observe that the higher numerical cost w.r.t. the bivariate experiments is due to the evaluation of the scaled Chebyshev-Vandermonde matrix at QMC points. {\textcolor{black}{In particular, the determination of QMC rule in the element $\Omega_3$ is more time consuming, since such is the in-domain routine on polyhedron w.r.t. that on union of balls. We can also see that for low degrees the compression time is negligible, while it should be taken into account for larger ADE}}.

We conclude by listing in Table {\ref{TAB_QMC_W}} the stability parameter $\|\mathbf{w}\|_1/|\sum_j w_j|\approx \|\mathbf{w}\|_1/vol(\Omega_s)$, $s=1,2$. The experiments show that these quantities are close to the optimal value 1, ensuring a good stability of the rule. Notice that in this case (\ref{stab}) is clearly an overestimate of the size of $\|\mathbf{w}\|_1$. This can be ascribed to the fact that the QMC moments of the product Chebyshev basis approximate the Lebesgue moments of the basis (at least when the number $L$ of QMC points is large) 
so that, rather, a bound close to (\ref{BM}) and consequently a substantial boundedness of the stability parameter is expected. 
\vskip0.5cm
\begin{Remark}
The possibility of integrating numerically on very complex  spatial elements without a difficult accurate tracking of the geometrical shape, makes in principle QMC an appealing approach, alternative to traditional cubature methods. Its Cheap version is able to reduce substantially the amount of computations, when a low-cost and stable integration formula is needed for several integrations on the same element. On the other hand, the applicability within FEM/VEM  polytopal discretizations is questionable, due to the relatively high cost of the initial construction of the QMC nodes and moments. Nevertheless, a possible interest may arise for example with parametric PDEs, where a given discretization can be used repeatedly changing the parameters, and certain preparatory computations can be made once per element and used in different simulations. 
   
\end{Remark}

\vskip0.5cm 
\noindent
{\bf Acknowledgements.} 

Work partially supported by the DOR funds of the University of Padova and by the INdAM-GNCS 2025 Project ``Polynomials, Splines and
Kernel Functions: from Numerical Approximation to Open-Source Software''. 
This research has been accomplished within the Community of Practice 
``Green Computing" of the Arqus European University Alliance, the RITA ``Research ITalian network on Approximation", and the SIMAI Activity Group ANA\&A.

\end{document}